\documentclass[11pt, leqno]{article}

\usepackage[utf8]{inputenc} 
\usepackage[a4paper,margin=1.2901in]{geometry} 
\usepackage{lmodern} 
\usepackage{microtype} 
\usepackage[english]{babel} 

\usepackage{amsmath, amsthm, amscd, amsfonts, amssymb} 
\usepackage{mathtools} 
\usepackage{mathrsfs} 
\usepackage{dsfont} 
\usepackage{latexsym} 
\usepackage{bm} 
\usepackage{accents} 
\usepackage{xfrac} 
\usepackage{slashed} 

\usepackage{graphicx} 
\usepackage{xcolor} 

\usepackage{enumitem} 
\usepackage{etoolbox} 
\usepackage{changepage} 

\usepackage{mdframed} 
\usepackage{needspace} 

\usepackage{amsrefs} 
\usepackage[bookmarksnumbered, colorlinks, plainpages]{hyperref} 
\usepackage{cleveref} 
\hypersetup{
    colorlinks=true,       
    linkcolor=red,         
    urlcolor=blue,         
    citecolor=green,       
    filecolor=magenta,     
    linkbordercolor=red,   
    urlbordercolor=blue,   
    citebordercolor=green, 
    filebordercolor=magenta
}

\newlength{\abstractmargin}
\usepackage{titlesec}
\titlespacing*{\section}{0pt}{12pt plus 4pt minus 2pt}{6pt plus 2pt minus 2pt}
\titlespacing*{\subsection}{0pt}{9pt plus 4pt minus 2pt}{6pt plus 2pt minus 2pt}

\makeatletter
\renewenvironment{proof}[1][\proofname]{\par
  \pushQED{\qed}%
  \normalfont \topsep6\p@\@plus6\p@\relax
  \trivlist
  \item[\hskip\labelsep
    \bfseries
    #1\@addpunct{.}]\ignorespaces
}{%
  \popQED\endtrivlist\@endpefalse
  \addvspace{6\p@\@plus6\p@\relax}
}
\makeatother

\newtheoremstyle{lmodern}
  {5pt}
  {5pt}
  {\itshape}
  {}
  {\bfseries}
  {.}
  {.5em}
  {\thmname{\footnotesize\MakeUppercase{#1}}\thmnumber{ \normalsize\textbf{#2}}\thmnote{ \normalsize\textbf{(#3)}}}

\theoremstyle{lmodern}
\newtheorem{theorem}{Theorem}[section]
\newtheorem{corollary}{Corollary}[theorem]
\newtheorem{lemma}[theorem]{Lemma}
\newtheorem{definition}{Definition}
\newtheorem{remark}{Remark}[section]
\newtheorem{proposition}[theorem]{Proposition}

\theoremstyle{remark}
\newtheorem{note}{\normalsize\textbf{Note}}[section] 

\numberwithin{equation}{section}

\newcommand{\norm}[1]{\left\Vert#1\right\Vert} 
\newcommand{\p}{\partial} 
\newcommand{\abs}[1]{\left\vert#1\right\vert} 
\newcommand{\Real}{\mathbb{R}} 
\newcommand{\ubi}{\textbf{\textit{u}}} 
\newcommand{\ybi}{\textbf{\textit{y}}} 
\newcommand{\Xd}{\textbf{\textit{X}}} 
\newcommand{\w}{\bm{\omega}} 
\newcommand{\al}{\bm{\alpha}} 
\newcommand{\spc}{\kern 0.08333em} 

\newcommand{\diagfrac}[2]{
    \scriptsize{\raise0.6ex\hbox{$#1$} \!\mathord{\left/{\vphantom{#2}}\right.\kern-\nulldelimiterspace}\!\lower0.6ex\hbox{$#2$}}%
}

\newcommand{\myframed}[3][0.4\baselineskip]{
    \vspace{#1} 
    \begin{mdframed}[
        backgroundcolor=yellow!5!white,
        linecolor=cyan!60!black,
        linewidth=0.5pt,
        roundcorner=10pt,
        skipabove=#2,
        skipbelow=3pt,
        innertopmargin=0pt,
        innerbottommargin=5pt,
        innerleftmargin=7pt, 
        innerrightmargin=8pt, 
        leftmargin=-8pt, 
        rightmargin=-8pt, 
    ]
    #3
    \end{mdframed}
}

\makeatletter
\renewcommand{\tableofcontents}{%
    \@starttoc{toc}%
}
\makeatother

\title{Uniform Boundedness \\ of Homogeneous Incompressible Flows in $\Real^3$}
\author{{Ulisse Iotti}\footnote{ulisse.iotti@globalscienceresearch.org \smallskip}}
\date{October 3, 2024}

\allowdisplaybreaks
\begin{document}

\maketitle

\begin{abstract}
    This paper investigates the extendability of local solutions for incompressible 3D Navier-Stokes and 3D Euler problems, with initial data $\ubi_0$ in the Sobolev space $H^s (\Real^3)$, where $s$ ensures the existence and uniqueness of classical solutions. A geometric decomposition of the configuration space, identified by the orthogonality between the solution $\ubi$ and the pressure forces $\nabla p$, splits the problem into two simpler subproblems, which enable the uniform boundedness of the solution in each component of the partition, thereby ensuring the extendability of the solution.
\end{abstract}

\vspace{1em}

\let\oldthefootnote\thefootnote 
\renewcommand{\thefootnote}{} 
\footnotetext{This work is distributed under the Creative Commons Attribution-NonCommercial-ShareAlike 4.0 International license for non-commercial use \href{https://creativecommons.org/licenses/by-nc-sa/4.0/?ref=chooser-v1}{(CC BY-NC-SA 4.0)}. This document is an unsigned copy. The original is signed in PADES-LT format, compliant with eIDAS, by Ulisse Iotti on 06/03/2025. \href{https://tinyurl.com/globalscienceresearch}{Verify the original signed document here}.}
\let\thefootnote\oldthefootnote 





\medskip
\medskip

\section*{Contents}
\begin{adjustwidth}{\abstractmargin}{\abstractmargin}
    \tableofcontents
\end{adjustwidth}

\clearpage
\section{Introduction}
In this paper we consider  the initial value problem for the homogeneous incompressible 3D Euler equations\
\begin{equation}\label{Euler}
\begin{cases}
  \partial_t \ubi + (\ubi \cdot \nabla)\ubi + \nabla p = 0 \qquad   \qquad& \, (x,t) \; \in \; \Real^3 \times [0,T),
  \\
  \nabla \cdot \ubi = 0 & \,  (x,t) \; \in \; \Real^3 \times [0,T),
  \\
  \ubi(x,0) = \ubi_0(x) & \, x \in \Real^3
\end{cases}
\end{equation}
%
and the homogeneous incompressible 3D Navier-Stokes equations\
\begin{equation}\label{NSE}
\begin{cases}
  \partial_t \ubi + (\ubi \cdot \nabla)\ubi - \mu \spc \Delta \ubi + \nabla p = 0 \qquad   \qquad& \, (x,t) \; \in \; \Real^3 \times [0,T),
  \\
  \nabla \cdot \ubi = 0 & \,  (x,t) \; \in \; \Real^3 \times [0,T),
  \\
  \ubi(x,0) = \ubi_0(x) & \, x \in \Real^3
\end{cases}
\end{equation}
%
Let $\ubi(x,t) = (u_1(x,t), u_2(x,t), u_3(x,t))$ denote the velocity field, and let $p(x,t)$ represent the pressure function of an ideal fluid. Here, $\mu$ is the kinematic viscosity, which is a positive constant.
\\
A salient feature of the Euler and Navier-Stokes equations is the inclusion of the pressure term, which is not governed by a dynamical equation. Indeed, pressure is not an independent variable; rather, it functions as a Lagrange multiplier that enforces the incompressibility condition $\nabla \cdot \ubi = 0$.
In problems \eqref{Euler} and \eqref{NSE}, the kinematic pressure $p$  is obtained by solving Poisson's problem, which is derived by taking the divergence of the evolution equations \eqref{Euler} or \eqref{NSE}:
    \begin{equation}\label{Pressure}
        \begin{cases}
            -\Delta p = \nabla\cdot[(\ubi \cdot \nabla)\ubi] = \nabla \cdot (\ubi \otimes \ubi)=  Tr (\nabla \ubi)^2  
            \\
            \mathop {\lim }\limits_{\abs{x} \to \infty} p = 0
        \end{cases}
    \end{equation}
In the whole space with decaying boundary conditions, the pressure is given by:\\
$$
    p=\sum\limits_{j,k} R_i R_j (u_i u_j)
$$
where \mbox{$R_i= \partial_i (-\Delta)^{-1}$} are Riesz transforms. The pressure is defined up to a time-dependent constant, and the usual choice involves zero average pressure. It can be noted that the pressure forces are not of a local nature.

\smallskip
The Navier-Stokes \eqref{NSE} and Euler \eqref{Euler} equations describe the behavior of incompressible flows and play a central role in the theory of nonlinear partial differential equations. Among the fundamental open questions is the possibility of globally extending local regular solutions in the case of finite energy, or determining whether blow-up can occur in finite time.\\
Numerous authors who have studied the problem in depth have identified several discriminating properties of the local solutions, the boundedness of which guarantees the extendability, while the unboundedness leads to the formation of singularities.

To obtain an existence theorem, regardless of the functional space utilized, an approximating sequence of the solution will always be constructed, whose terms generally exhibit global existence in time; however, this global existence is lost in the transition to the limit when obtaining the final solution. This occurs because it is necessary to traverse a \emph{bottleneck} that cannot be crossed while preserving the approximations. This constraint is associated with some form of compactness or dominated convergence, which is controlled by the nonlinear term.\

It is hypothesized that this \emph{bottleneck} does not represent an inevitable impediment, but rather a technical limitation stemming from the conventional treatment of the nonlinear term as a straightforward product of functions in norm estimation. This approach inevitably neglects the Helmholtz projector and the intrinsic structure of the nonlinear term, which, if \emph{untied} like a trivial knot, could significantly simplify the problem.

The present study addresses the extendability of a family of classical local solutions for finite energy, where sufficient regularity of initial data permits simplified notation.
No explicit criterion for convergence or blow-up dynamics is proposed. The adopted approach involves a geometric decomposition of the problem into two analogous yet elementary subproblems, defined through an appropriate two-component partition of the configuration space.

\medskip
\medskip
\medskip
\section{The main result}
The main results of the paper are the following:

\myframed[0.4\baselineskip]{1pt}{
\begin{theorem}\label{Main01}
Let $\ubi(x,t)$ be a solution to \eqref{Euler} with initial condition $\ubi_0 \in H^s(\Real^3)$ where $s > \diagfrac{3}{2} + 2$. Then,
    \begin{equation}
    \norm{\ubi(\cdot,t)}_{\infty} = \norm{\ubi_0}_{\infty} \quad \forall t \in [0,T^{\ast})
    \end{equation}
    Hence, the solution $\ubi$ can be extended for all $t \geq 0$.
\end{theorem}
}

\myframed[0.4\baselineskip]{1pt}{
\begin{theorem}\label{Main02}
Let $\ubi(x,t)$ be a solution to \eqref{NSE} with initial condition $\ubi_0 \in H^s(\Real^3)$ where $s > \diagfrac{3}{2}-1$. Then,
    \begin{equation}
    \norm{\ubi(\cdot,t)}_{\infty} \leq \norm{\ubi(\cdot,\varepsilon)}_{\infty} \quad \forall t \in [\varepsilon,T^{\ast})
    \end{equation}
    Hence, the solution $\ubi$ can be extended for all $t \geq 0$.
\end{theorem}
}

\begin{remark}
  The aim of this work is not to determine the minimal regularity necessary, but to adopt the lowest level of regularity sufficient to optimize notational simplicity.
\end{remark}

\medskip
\subsection{Ideas of the proof}
To investigate the scenario of infinite extendability of solutions, we consider the velocity field $\ubi(x,t)$ as a solution, for instance, to the Euler problem \eqref{Euler}.
        $$D_t \ubi =\p_t \ubi + (\ubi\cdot\nabla) \ubi$$
is the material derivative of $\ubi$, expressing the variation of $\ubi$ over time along the stream lines. We can rewrite \eqref{Euler} as
        $$D_t \ubi = -\nabla  p$$
and by scalar multiplication with $\ubi$, we obtain
        $$\frac{1}{2} D_t \abs{\ubi \spc }^2 = -\nabla p \cdot \ubi$$
which expresses the variation of energy along the stream lines.
\\
If $\nabla p$ were always orthogonal to $\ubi$, we would have $D_t \abs{\ubi \spc }^2 \equiv 0$, leading to energy conservation along all stream lines and $\norm{\ubi(t)}_{\infty} \equiv \norm{\ubi_0}_{\infty}$.
\\
The condition $\nabla p \perp \ubi$ for all $(x,t) \in \Real^3 \times [0,T^*)$ turns out to be overly restrictive for general solutions $\ubi$.
Nevertheless, we impose this constraint by an auxiliary PDE with a new unknown field, thus introducing additional degrees of freedom into the constraint expression.
In this PDE, $\ubi$ and its spatial derivatives will act as coefficients.
The domain of existence of the solution to the auxiliary equation induces a two-component partition of $\Real^3$, which is preserved under the \textit{flow map} diffeomorphism.
This allows the problem to be split into two simpler problems and enables a straightforward characterization of $\ubi$ in each component of the partition of $\Real^3 \times [0, T^*)$.
\\

%
%
\medskip
\medskip
\section{Preliminaries and known results}
\subsection{Local existence and uniqueness theorems in the Eulerian approach}

In this work, we state only the classical existence and uniqueness theorems for Sobolev spaces $H^s(\Real^n)$, as they suffice for our analysis and no results in more general functional spaces are required.
For each theorem, we will provide bibliographic references to explore the most recent developments in much more general functional spaces.

\myframed[0.4\baselineskip]{1pt}{
\begin{theorem}[\small{Theorem of local existence of classical solutions for incompressible Navier-Stokes problems in $H^s(\Real^n)$ for $s > \diagfrac{n}{2} - 1$}] \label{T:KATO_NSE}
\hfill
\begin{enumerate}[label=(\roman*)]
    \item Let  $s > \diagfrac{n}{2} - 1$. For $\ubi_0 \in H^s(\mathbb{R}^n)$ with $\nabla \cdot \ubi_0 = 0$, there exist a positive $T^*$ and a (unique) weak solution $\ubi$ for the incompressible Navier-Stokes  problem \eqref{NSE} on $(0, T^*) \times \Real^n$ such that $\ubi(0, x) = \ubi_0(x)$. Moreover, the solution $\ubi$ satisfies the following:
    \begin{enumerate}[label=(\alph*)]
        \item $\ubi \in C\left((0, T^*), H^s\Real^n \right)$
        \item $\ubi \in C^\infty \left((0, T^*), H^{\infty}(\Real^n)\right)$
        \item $\ubi \in C^\infty \left((0, T^*) \times \Real^n\right)$
        \item $\ubi \in L^2\left((0, T), H^{s+1}(\Real^n)\right)$ for all \, $T \in (0, T^*)$
    \end{enumerate}

    \item The maximal existence time $T^*$ of the smooth solution $\ubi$ is finite if and only if\\ $\lim_{t \to T^*} \norm{\ubi(t)}_\infty = \infty$.
\end{enumerate}
\end{theorem}
}

In local existence theorems for the Navier-Stokes equations \eqref{NSE} and the Euler equations \eqref{Euler}, we typically analyze the projected equation using the Helmholtz projector, which maps onto divergence-free vector fields and eliminates the pressure gradient.

For the Navier-Stokes equations, the integral form of the equation and the semigroup generated by the projected Laplacian (the Stokes operator) play a crucial role, as the latter provides regularization properties essential for this formulation.

Existence theorems, in general, depend on approximation or regularization strategies, leading to the construction of a sequence that converges, in some topology, to the solution. This sequence is frequently generated --- i.e., its existence is proven --- using techniques related to some form of fixed-point argument. Typically, passing to the limit requires a form of compactness, obtained through a priori estimates or uniform properties tied to the nonlinear term. However, ensuring this compactness typically necessitates sacrificing certain properties of the sequence elements, such as global existence or some degree of regularity.

Theorems in Sobolev spaces are essentially derived from the works of \cite{Kato1964}, \cite{Kato1972}, \cite{Kato1984}, \cite{Kato1988} and \cite{Kato1992}. Kato's mild solutions paved the way for developments in Besov spaces, as seen in \cite{Cannone1995}, \cite{Cannone1997}, \cite{Bahouri2011}  and later extended to $BMO^{-1}$ spaces, pioneered by \cite{Koch2001}.

For a comprehensive overview of existence and uniqueness theorems for the incompressible Navier-Stokes equations in Sobolev, Besov, homogeneous Besov, $BMO$, and $BMO^{-1}$ spaces, see \cite{Lemarie2002}, \cite{Kozono2019}. 

\myframed[0.4\baselineskip]{1pt}{
\begin{theorem}[\small{Theorem of local existence of classical solutions for incompressible Euler problems in $H^s(\Real^n)$ for $s > \diagfrac{n}{2} + 1$}] \label{T:LocalEuler}
\hfill
\begin{enumerate}[label=(\roman*)]
    \item Let $k \geq 0$ and $s > \diagfrac{n}{2} + k + 1$. For $\ubi_0 \in H^s(\mathbb{R}^n)$ with $\nabla \cdot \ubi_0 = 0$, there exist a positive $T^*$ and a classical solution $\ubi$ for the incompressible Euler problem \eqref{Euler} on $(-T^*, T^*) \times \Real^n$ such that $\ubi(0, x) = \ubi_0(x)$ ad $0 \leq r \leq k$. Moreover, the solution $\ubi$ satisfies the following:
    \begin{enumerate}[label=(\alph*)]
        \item $\ubi \in C^r\left((-T^*, T^*), H^{(s-r)}(\Real^n)\right)$
        \item $\ubi \in C^r\left((-T^*, T^*), C^{(k+1-r)}(\Real^n)\right)$
    \end{enumerate}

    \item The maximal existence time $T^*$ of the smooth solution $\ubi$ is finite if and only if\\ $\lim_{t \to T^*} \norm{\ubi(t)}_\infty = \infty$.
\end{enumerate}
\end{theorem}
}

For a comprehensive overview of existence and uniqueness theorems for the incompressible Euler equations in Sobolev, Besov, homogeneous Besov, Besov-Herz, and Hölder spaces, see \cite{Bardos2007} and \cite{Lucas2017}.

\medskip
\subsection{Sufficient conditions for extendability}

A common feature of existence theorems is the unavoidable dependence on the time integral of some norm of $\nabla \ubi$ or, $\w$ whether one uses a form of the logarithmic Sobolev inequality or studies the vorticity equation directly.
The behavior of these quantities provides sufficient criteria for extendability or necessary conditions for blow-up.
From Theorems \ref{T:KATO_NSE} and \ref{T:LocalEuler}, the following theorem follows directly.

\myframed[0.4\baselineskip]{1pt}{
\begin{theorem}[\small{$L^\infty$ gradient control and global existence}] \label{T:BMK_gradient}
\hfill\par\noindent
Let $\ubi_0 \in H^s(\Real^n)$ for some $s > \diagfrac{n}{2} + 1$, so that there exists a classical solution\\
$\ubi \in C\left([0, T^*), H^s(\Real^n)\right)$ to the incompressible Euler or Navier-Stokes problem.\\
Then:
\begin{enumerate}[label=(\roman*)]
    \item If for any $T > 0$ there exists $0 < M_1 < \infty$ such that
    \[
    \int_0^T \norm{\nabla \ubi(t)}_\infty \, dt \leq M_1,
    \]
    then the solution $\ubi$ exists globally in time and $\ubi \in C\left([0, \infty), H^s(\Real^n)\right)$.

    \item If the maximal time $T^*$ of the existence of the solution $\ubi$ is finite, then necessarily
    \[
    \lim_{t \to T^*} \int_0^T \norm{\nabla \ubi(t)}_\infty \, dt = \infty.
    \]
\end{enumerate}
\end{theorem}
}

\smallskip
This is followed by the famous theorem of BKM \cite{Beale1984}.
\smallskip

\myframed[0.4\baselineskip]{1pt}{
\begin{theorem}[\small{BKM. $L^\infty$ vorticity control and global existence}]\label{BKM}
\hfill\par\noindent
Let $\ubi_0 \in H^s(\Real^3)$ for some $s \geq \diagfrac{3}{2} + 1$, so that there exists a classical solution\\
$\ubi \in C\left([0, T^*), H^s(\Real^3)\right)$ to the incompressible Euler or Navier-Stokes problem.\\
Let $\w = \nabla \times \ubi$, then:
\begin{enumerate}[label=(\roman*)]
    \item If for any $T > 0$ there exists $0 < M_1 < \infty$ such that
    \[
    \int_0^T \norm{\,\w(t)}_\infty \, dt \leq M_1,
    \]
    then the solution $\ubi$ exists globally in time and $\ubi \in C\left([0, \infty), H^s(\Real^3)\right)$.
    \item If the maximal time $T^*$ of the existence of the solution $\ubi$ is finite, then necessarily
    \[
    \lim_{t \to T^*} \int_0^T \norm{\,\w(t)}_\infty \, dt = \infty.
    \]
\end{enumerate}
\end{theorem}
}

Vorticity is a typical concept defined in $\Real^3$ can be extended to higher dimensions using the exterior derivative of the covelocity 1-form $\ubi^*$. The vorticity 2-form $\w$ is given by:
$$
\w := d\ubi^*
$$
where $\w_{i,j} = \nabla_i \ubi_j - \nabla_j \ubi_i$. Alternatively, vorticity can be viewed as the antisymmetric part of the velocity gradient tensor $\nabla \ubi$:
$$
\w = (\nabla \ubi - \nabla \ubi^T).
$$
A significant version of \ref{BKM} was obtained by \cite{Kozono2000a} in $L^p$ spaces for $p>1$ and \cite{Kozono2000}, \cite{Kozono2002} proved that the $\norm{\cdot}_\infty$ norm can be replaced by the norm in the $BMO$ space. This generalization is interesting because some crucial Sobolev embedding theorems can be applied to the $BMO$ space, but not to the $L^\infty$ space. 

\smallskip
The BKM non-blow-up criterion has evolved in various forms, often by incorporating geometric aspects of the solution. It was first introduced in the work of \cite{Constantin1993}, where the direction of vorticity was analyzed. This approach was further developed in \cite{Constantin1996} and later extended in different directions by many others, such as \cite{Gibbon1997}, \cite{Cordoba2001a}, \cite{Deng2005}, \cite{Berselli2009}, \cite{Beirao2016} and \cite{Chae2016}. 
\

\smallskip
For recent developments on quantitative bounds and critical behavior in solutions to the incompressible Navier-Stokes equations, see \cite{Gallagher2018}, \cite{Tao2021} and \cite{Barker2022}.

\subsection{Representation by the Lamb Vector}\label{note01}
    By decomposing $\nabla \ubi$ into its symmetric and antisymmetric parts, the equation \eqref{NSE} or \eqref{Euler}, for $\mu \geq 0$, can be rewritten in the form
    \begin{equation}\label{Ortho52}
        \partial_t \ubi + \frac{1}{2} \nabla \abs{\ubi \spc }^2 - \ubi \times \w  - \mu\Delta \ubi + \nabla p = 0.
    \end{equation}

\begin{note}[\small\textbf{Helmholtz-Hodge decomposition of the Lamb vector $\ubi \times \w$}]
    Taking the divergence of (\ref{Ortho52}) yields
    \begin{equation*}
        \nabla \cdot (\ubi \times \w) = \frac{1}{2} \Delta \abs{\ubi \spc }^2 + \Delta p.
    \end{equation*}
    Consequently, we obtain
    \begin{equation*}
        \ubi \times \w =  \frac{1}{2} \nabla \abs{\ubi \spc }^2 + \nabla p + {\textbf{\textit{a}}}.
    \end{equation*}
    where, ${\textbf{\textit{a}}}$ denotes a vector field with $\nabla \cdot {\textbf{\textit{a}}} = 0$ \, and \, $\nabla \times {\textbf{\textit{a}}} = \nabla \times (\ubi \times \w)$.\\
\end{note}

\medskip
\subsection{Foundations of the Lagrangian approach for incompressible fluids}
\label{subsec:Lagrangian}

The Lagrangian framework provides a geometric perspective for analyzing incompressible fluid dynamics by tracking the motion of individual fluid particles. This approach intrinsically incorporates the divergence-free constraint and reveals deep connections with differential geometry and dynamical systems.

\subsubsection*{Particle Trajectories and Well-Posedness}
Let $\ubi: \Real^n \times [0,T) \to \Real^n$ be a velocity field solving the Euler \eqref{Euler} or Navier-Stokes \eqref{NSE} equations. The fundamental kinematic law is encoded in the system of ordinary differential equations:
\begin{equation}\label{eq:particle_ode}
    \frac{\partial}{\partial t} \Xd(\al, t) = \ubi(\Xd(\al, t), t), \quad \Xd(\al, 0) = \al,
\end{equation}
where $\Xd(\al, t)$ denotes the position at time $t$ of the fluid particle initially at $\al \in \Real^n$. Under the regularity assumption $\ubi \in C^r(\Real^n \times [0, T^*))$ with $r \geq 1$, the Cauchy-Lipschitz theorem guarantees:
existence and uniqueness of $\Xd \in C^r(\Real^n \times [0, T^*))$, continuous dependence on initial conditions, preservation of topological structure (no particle crossings).\\
This establishes $\Xd$ as a \textit{flow map} generating particle trajectories.

\subsubsection*{Incompressibility as Volume Preservation}
The divergence-free condition $\nabla \cdot \ubi = 0$ translates into a geometric constraint on the flow map. Let $J(\al, t) := \det(\nabla_{\al} \Xd(\al, t))$ be the Jacobian determinant of the deformation gradient tensor. Then:
\begin{equation}\label{eq:jacobian_constraint}
  J(\al, t) = 1 \quad \forall t \geq 0, \; \forall \al \in \Real^n,
\end{equation}
implying $\Xd(\cdot, t)$ is a \textit{volume-preserving diffeomorphism}. Physically, this means:
fluid parcels evolve without compression/expansion, the pushforward map $\Xd(\cdot,t)_*$ preserves Lebesgue measure, inverse map $\Xd^{-1}(\cdot,t)$ exists and is smooth.

\subsubsection*{Material Derivatives and Acceleration}
The \textit{material derivative} quantifies time evolution along particle paths. For any Eulerian function $f(x,t)$, its Lagrangian rate of change is:
\begin{equation}\label{eq:material_derivative}
  D_t f(x, t) := \left.\frac{d}{dt}\right|_{\al \text{ fixed}} f(\Xd(\al,t),t) = \partial_t f + \ubi \cdot \nabla f,
\end{equation}
where $x = \Xd(\al,t)$. For the velocity field itself, this gives:
\begin{equation}\label{eq:acceleration}
  D_t \ubi = \partial_t \ubi + \ubi \cdot \nabla \ubi = \frac{\partial^2 \Xd}{\partial t^2},
\end{equation}
directly linking to the momentum equations. The equivalence in \eqref{eq:acceleration} arises from differentiating $\ubi(\Xd(\al,t),t)$ via the chain rule.

\subsubsection*{Vorticity Transport and Lagrangian Coherence}
A fundamental result is the \textit{vorticity transport formula}, describing the evolution of $\w = \nabla \times \ubi$:
\begin{equation}\label{eq:vorticity_transport}
    \w(\Xd(\al,t),t) = \nabla_{\al} \Xd(\al,t) \, \w_0(\al),
\end{equation}
where $\w_0$ is the initial vorticity. This reveals that vorticity vectors are \textit{Lie-transported} along particle paths.


For \eqref{Euler}, existence and uniqueness theorems can be directly established for the the system \eqref{eq:particle_ode} formulated as integro-differential equations in H\"{o}lder spaces, using potential theory \cite{Majda2002}.\\
Highly fertile theoretical frameworks are those based on Lagrangian and Hamiltonian mechanics, leveraging the Frobenius theorem or analyzing structures on diffeomorphism groups \cite{Arnold1966}, \cite{Serre1984}. However, within these approaches, it is generally challenging to directly establish existence theorems except in particular cases with special symmetries \cite{Tao2016}, \cite{Tao2018}, \cite{Tao2020}, \cite{Shimizu2023}, among many others who have investigated the potential for finite-time blow-up.

\medskip
\medskip
\section{Partition of $\Real^3 \times [0, T^*)$ leading to problem simplification}\label{SezSet}

\medskip
\begin{definition}\label{zeriOmega}
    Let $\ubi \in C\left([0, T^*), C^2(\Real^n)\right)$ solution to \eqref{NSE} or \eqref{Euler}, and let \,
    \(\w = \nabla \times \ubi\).\\
    For $t \in [0, T^*)$ we define
    \begin{equation}\label{zeroOm}
        \Omega_t = \big\{ x \in \Real^3 \mid \, \w(x, t) = 0 \big\}.
    \end{equation}
    Since \(\w(\cdot, t) \in C^1(\Real^3)\), it follows that \( \Omega_t \) is a \( C^1 \)-regular closed set.
    \medskip
    \begin{remark}\label{zeriOmega02}
        By \eqref{eq:vorticity_transport}, the sets $\Omega_t$ are mutually diffeomorphic for different values of $t$ via the map $\Xd(\cdot, t)$,
        and by virtue of \eqref{eq:jacobian_constraint}, they share the same measure. Particularly, we can write $\Omega_t = \Xd(\Omega_0, t)$.\\
        $\Omega_t \times [0,T^*) = \bigcup_{x \in \Omega_0} \Xd(x, [0,T^*))$ is a stream tube.
    \end{remark}
    \medskip
    \begin{remark}\label{zeriOmega03}
        Similarly, the sets $\Real^3 \setminus \Omega_t$, $\overline{\Real^3 \setminus \Omega_t}$ and $\partial \spc \Omega_t$ are mutually diffeomorphic via $\Xd(\cdot, t)$, share the same measure and identify stream tubes.
    \end{remark}
    \medskip
\end{definition}

\subsection{The $(\ubi, \nabla p)$\spc-\spc orthogonality set}

\myframed[0.4\baselineskip]{1pt}{
\begin{lemma}\label{Ortho01}
    Let $\ubi \in C([0, T^*), H^s(\Real^3))$ for $s > \diagfrac{3}{2}+2$, solution to \eqref{NSE} or \eqref{Euler} and $t \in [0, T^*)$. Then
    \begin{equation}
        \ubi \perp \nabla p \quad  \text{ in } \quad \Real^3 \setminus \Omega_{t}.
    \end{equation}
\end{lemma}
}

\begin{proof}
    Let $t \in [0, T^*)$. To ensure that $\ubi \perp \nabla p$, we set
    \begin{equation}\label{Ortho03}
        \nabla p = \ubi \times \nabla \theta,
    \end{equation}
    where $\nabla \theta$ is a vector field to be determined. Then, $\ubi \perp \nabla p$ holds wherever $\nabla \theta$ exists.\\
    By vector identity
    \begin{equation*}
        \nabla \times (\varphi \spc \ybi) = \varphi \spc (\nabla \times \ybi) - \ybi \times \nabla \varphi
    \end{equation*}
    equation~\eqref{Ortho03}, can be rewritten as
    \begin{equation}\label{Ortho0301}
        \nabla p = \ubi \times \nabla \theta = \theta \spc \w - \nabla \times \theta \ubi.
    \end{equation}
    Taking the pointwise scalar product of both sides of~\eqref{Ortho0301} with $\w$, we obtain
    \begin{equation}\label{Ortho0302}
        \w \cdot \nabla p = \theta \spc {\abs{\spc \w}}^2 - \w \cdot (\nabla \times \theta \ubi).
    \end{equation}
    Rearranging terms yields
    \begin{equation}\label{Ortho0303}
        \theta = \frac{\w \cdot \nabla p}{{\abs{\spc \w}}^2} + \frac{\w \cdot (\nabla \times \theta \ubi)}{{\abs{\spc \w}}^2}.
    \end{equation}
    From this expression, it follows that the existence of $\theta$ is not guaranteed in $\Omega_t$, since ${\abs{\spc \w}}^2 = 0$ there.
    \\

    Next, consider $\Real^3 \setminus \Omega_t$, where $\w \neq 0$. Taking the divergence of~\eqref{Ortho0301}, we find
    \begin{equation}\label{Ortho0304}
        \Delta p = \nabla \theta \cdot \w.
    \end{equation}
    This is a non homogeneous first-order linear partial differential equation for $\theta$. Since $\w \in C^1(\Real^3)$, $\Delta p \in C^1(\Real^3)$, and $\w \neq 0$ in $\Real^3 \setminus \Omega_t$, with both $\w$ and $\Delta p$ vanishing at infinity and $\w=0$ in $\partial \spc \Omega_t$ and $\Omega_t$ is a $C^1$-regular set, the method of characteristics guarantees the existence of a solution $\theta \in C^1(\Real^3 \setminus \Omega_t)$.\\
    Consequently, for all $t \in [0,T^*)$, we have
    \begin{equation}
        \nabla p = \ubi \times \nabla \theta \quad \text{in } \Real^3 \setminus \Omega_t,
    \end{equation}
    and by construction,
    \begin{equation}\label{Ortho04}
        \ubi \cdot \nabla p = 0 \quad \forall x \in \Real^3 \setminus \Omega_t.
    \end{equation}
\end{proof}

\medskip
\myframed[0.4\baselineskip]{1pt}{
\begin{corollary}\label{Frontier}
    $\nabla p \cdot \ubi = 0 \quad \forall x \in \partial \spc \Omega_t$.
\end{corollary}
}

\begin{proof}
    Since $(\nabla p \cdot \ubi) \in C(\Real^3 \times [0, T^*))$, continuity yields
    \begin{equation*}
        \nabla p \cdot \ubi = 0 \quad \text{in} \quad \Real^3 \setminus \Omega_t \quad \Rightarrow \quad \nabla p \cdot \ubi = 0 \quad \text{in} \quad \overline{\Real^3 \setminus \Omega_t}.
    \end{equation*}
\end{proof}

\medskip
\myframed[0.4\baselineskip]{1pt}{
\begin{corollary}\label{Frontier02}
    If $\mathring{\Omega}_t = \emptyset$, \, then \, $\nabla p \cdot \ubi = 0 \quad \forall x \in \Real^3$.
\end{corollary}
}

\begin{proof}
    This is a direct consequence of Lemma \ref{Ortho01} and Corollary \ref{Frontier}.
\end{proof}

\medskip
\subsection{The zero-vorticity set}

\myframed[0.4\baselineskip]{1pt}{
\begin{proposition}\label{Vv01}
    Let $\ubi$ be as defined in Lemma \ref{zeriOmega}.
    In \mbox {$\Omega_t$}, we have
    \begin{equation*}
        \sup_{x \in \Omega_t } \abs{\ubi \spc }^2 \leq \sup_{x \in \partial \spc \Omega_t } \abs{\ubi \spc }^2
    \end{equation*}
\end{proposition}
}

\begin{proof}
    In $\Omega_t$, we have $\nabla \times \ubi = 0$, hence there exists a scalar potential $\phi$ such that $\ubi = \nabla \phi$. \\
    Since $\nabla \cdot \ubi = 0$, \spc $\phi$ is harmonic, and $(\nabla \phi)^2$ is subharmonic (see note \ref{note04}). \\
    For $\abs{\ubi \spc }^2 = (\nabla \phi)^2$ continuous and subharmonic in $\Omega_t$, the weak maximum principle holds:
    \begin{equation}\label{Harmonic01}
        \sup_{x \in \Omega_t} \abs{\ubi \spc}^2 \leq \sup_{x \in \partial \spc \Omega_t} \abs{\ubi \spc}^2.
    \end{equation}
\end{proof}

\medskip
\myframed[0.4\baselineskip]{1pt}{
\begin{corollary}\label{Vv02}
    We have \, $\norm{\ubi(\cdot,t)}_{\infty, \Omega_t} \leq \norm{\ubi(\cdot,t)}_{\infty, \overline{\Real^3 \setminus \Omega_t}}$ \; for \; $t \in [0, T^*)$
\end{corollary}
}

\begin{proof}
    This is a direct consequence of Proposition \ref{Vv01}.\\
\end{proof}

\begin{note}[\small\textbf{Subharmonicity of $\abs{\ubi \spc }^2$}]\label{note04}
\hfill\break
$\phi$ is harmonic in $\Omega_t$
$$
    \abs{\ubi \spc }^2 = (\nabla \phi)^2 = \sum_{i} \left(\frac{\partial \phi}{\partial x_i}\right)^2
$$
$$
    \Delta \abs{\ubi \spc }^2  = \sum_{j} \frac{\partial^2 \ubi^2}{ \partial x_j^2}= 2 \sum_{j} \left(\frac{\partial \ubi}{ \partial x_j}\right)^2 + 2 \ubi \cdot \sum_{j} \frac{\partial^2 \ubi}{ \partial x_j^2}
    = 2 \sum_{j} \left(\frac{\partial \ubi}{ \partial x_j}\right)^2 + \ubi \cdot \Delta \ubi
$$
$$
    \ubi \cdot \Delta \ubi = 2 \sum_{i} \frac{\partial \phi}{\partial x_i} \Delta \frac{\partial \phi}{\partial x_i} = 2 \sum_{i} \frac{\partial \phi}{\partial x_i} \frac{\partial \Delta \phi}{\partial x_i}= 0
$$
Thus
$$
     \Delta \abs{\ubi \spc }^2  = 2 \sum_{j} \left(\frac{\partial \ubi}{ \partial x_j}\right)^2 \geq 0.
$$
\end{note}

\clearpage
\section{Uniform boundedness of solutions for the Euler problem}
\myframed[0.4\baselineskip]{1pt}{
\begin{theorem}[Uniform Boundedness of the Euler Solution in $\Real^3$]\label{thm:BoundEuler}
  \hfill\break
  \begin{minipage}[t]{\linewidth}
    Let \(\ubi(\cdot, t)\) be a solution of problem \eqref{Euler} in \(H^s(\Real^3)\) with \(s > \diagfrac{3}{2}+2\). Then the following hold:
        \begin{enumerate}[label=\Alph*), labelsep=0.5cm, leftmargin=1.5cm, itemsep=1em]
            \item \(D_t \abs{\ubi \spc }^2 = 0\) \quad \text{for all} \quad \((x, t) \in \overline{\Real^3 \setminus \Omega_t} \times [0, T^*)\).
            \item $\norm{\ubi(\cdot,t)}_{\infty} = \norm{\ubi(\cdot,t)}_{\infty, \overline{\Real^3 \setminus \Omega_t}}  \quad \text{for all} \quad t \in [0,T^{\ast})$.
            \item $\norm{\ubi(\cdot,t)}_{\infty} = \norm{\ubi_0}_{\infty} \quad \text{for all} \quad t \in [0,T^{\ast})$
        \end{enumerate}
  \end{minipage}
  \\
  \\
  \hfill\break
  Hence, the solution $\ubi$ can be extended for all $t \geq 0$.
\end{theorem}
}

\begin{proof}[Proof( of Main Theorem \ref{Main01})]
    \hfill\\
    For $s > \diagfrac{3}{2}+2$, Theorem \ref{T:LocalEuler} ensures the existence of the solution $\ubi \in C([0,T^*), C^2(\Real^3) \cap H^s(\Real^3))$.
    \\
    \\
    (\textit{A}) We start with \eqref{Ortho52} and take the pointwise scalar product with \(\ubi\), yielding
    \begin{equation}\label{BoundEuler01}
       \frac{1}{2} \partial_t \abs{\ubi \spc }^2 = -\frac{1}{2} \nabla \abs{\ubi \spc }^2 \cdot \ubi + (\ubi \times \w) \cdot \ubi - \nabla p \cdot \ubi .
    \end{equation}
    By Lemma \ref{Ortho01}, Corollary \ref{Frontier} and Remark \ref{zeriOmega03}, we have \(\nabla p \cdot \ubi = 0\) for all \(x \in \overline{\Real^3 \setminus \Omega_t}  \times [0, T^*) \).\\
        Consequently,
    $$
        D_t \abs{\ubi \spc }^2 = 0 \quad \quad \forall (x,t) \, \in \, \overline{\Real^3 \setminus \Omega_t} \times [0, T^*).
    $$
    \\
    \\
    (\textit{B})
    \begin{equation*}
        \norm{\ubi(\cdot,t)}_{\infty} = \max \left\{ \norm{\ubi(\cdot,t)}_{\infty, \Omega_t}, \norm{\ubi(\cdot,t)}_{\infty, \overline{\Real^3 \setminus \Omega_t}} \right\} \quad \forall t \in [0, T^*)
    \end{equation*}
    By Corollary \ref{Vv02}, we obtain immediately
    \begin{equation*}
        \norm{\ubi(\cdot,t)}_{\infty} = \norm{\ubi(\cdot,t)}_{\infty, \overline{\Real^3 \setminus \Omega_t}} \quad \forall t \in [0, T^*).
    \end{equation*}
    \\
    \\
    (\textit{C}) By point (\textit{A}), energy is constant along all streamlines in \mbox {$\overline{\Real^3 \setminus \Omega_t} \times [0, T^*)$}. This leads to
    \begin{equation}\label{norm02}
        \norm{\ubi(\cdot,t)}_{\infty, \overline{\Real^3 \setminus \Omega_t}} = \norm{\ubi_0}_{\infty, \overline{\Real^3 \setminus \Omega_0}} \quad \forall t \in [0, T^*).
    \end{equation}
    Moreover, (\textit{B}) ensures that
    \begin{equation}\label{norm03}
        \norm{\ubi(\cdot,t)}_{\infty} = \norm{\ubi_0}_{\infty}
    \end{equation}
    Since $\lim_{t \to T^*} \norm{\ubi(t)}_\infty < \infty$, Theorem \ref{T:LocalEuler} \textit{(ii)} guarantees that the solution $\ubi$ can be extended for all $t \geq 0$.
    \\
\end{proof}

\clearpage
\section{Uniform boundedness of solutions for the Navier-Stokes problem}
\myframed[0.4\baselineskip]{1pt}{
\begin{theorem}[Uniform Boundedness of the Navier-Stokes Solution in $\Real^3$]\label{BoundNSE}
  \hfill\break
  \begin{minipage}[t]{\linewidth}
    Let $\ubi(\cdot, t)$ be a solution of problem \eqref{NSE} in $H^s(\Real^3)$ with $s > \diagfrac{3}{2}-1$ and $0<\varepsilon <T^*$. Then the following hold:
    \begin{enumerate}[label=\Alph*), labelsep=0.5cm, leftmargin=1.5cm, itemsep=1em]
            \item \(D_t \abs{\ubi \spc }^2 \leq  \mu \Delta \abs{\ubi \spc }^2 \) \quad \text{for all} \quad \((x, t) \in \overline{\Real^3 \setminus \Omega_t} \times [\varepsilon, T^*)\).
            \item $\norm{\ubi(\cdot, t)}_{\infty} = \norm{\ubi(\cdot, t)}_{\infty, \overline{\Real^3 \setminus \Omega_t}}  \quad \text{for all} \quad t \in (0, T^{\ast})$.
            \item $\norm{\ubi(\cdot,t)}_{\infty} \leq \norm{\ubi(\cdot,\varepsilon)}_{\infty} \quad \forall t \quad \in [\varepsilon,T^{\ast})$
    \end{enumerate}
  \end{minipage}
  \\
  \\
  \hfill\break
  Hence, the solution $\ubi$ can be extended for all $t \geq 0$.
\end{theorem}
}

\begin{proof}[Proof of Main Theorem \ref{Main02}]
    \hfill\\
    For $s > \diagfrac{3}{2}-1$, Theorem \ref{T:KATO_NSE} ensures the existence of the solution $\ubi \in C^\infty((0,T^*), H^\infty(\Real^3))$.
    \\
    \\
    (\textit{A})  We start with \eqref{Ortho52} and take the pointwise scalar product with \(\ubi\), yielding
    \begin{equation}\label{BoundNSE01}
       \frac{1}{2} \partial_t \abs{\ubi \spc }^2 = -\frac{1}{2} \nabla \abs{\ubi \spc }^2 \cdot \ubi  + (\ubi \times \w) \cdot \ubi - \nabla p \cdot \ubi + \mu \Delta \ubi \cdot \ubi.
    \end{equation}
    By Lemma \ref{Ortho01}, Corollary \ref{Frontier} and Remark \ref{zeriOmega03}, we have \(\nabla p \cdot \ubi = 0\) for all \(x \in \overline{\Real^3 \setminus \Omega_t}  \times (0, T^*) \).\\
    Equation \eqref{BoundNSE01}, in $\overline{\Real^3 \setminus \Omega_t}  \times (0, T^*)$ can be rewritten as
    \begin{equation}\label{BoundNSE02}
       \frac{1}{2} D_t \abs{\ubi \spc }^2 =  \mu \Delta \ubi \cdot \ubi = \mu \frac{1}{2} \Delta \abs{\ubi \spc }^2 - \mu\abs{\nabla \ubi}^2,
    \end{equation}
    and thus
    \begin{equation}\label{BoundNSE03}
        D_t \abs{\ubi \spc }^2 \leq  \mu \Delta \abs{\ubi \spc }^2 \quad \quad \forall (x,t) \, \in \, \overline{\Real^3 \setminus \Omega_t} \times (0, T^*).
    \end{equation}
    \\
    \\
    (\textit{B})
    \begin{equation*}
        \norm{\ubi(\cdot, t)}_{\infty} = \max \left\{ \norm{\ubi(\cdot,t)}_{\infty, \Omega_t}, \norm{\ubi(\cdot,t)}_{\infty, \overline{\Real^3 \setminus \Omega_t}} \right\} \quad \forall t \in [0, T^*)
    \end{equation*}
    By Corollary \ref{Vv02}, we obtain immediately
    \begin{equation*}
        \norm{\ubi(\cdot,t)}_{\infty} = \norm{\ubi(\cdot,t)}_{\infty, \overline{\Real^3 \setminus \Omega_t}} \quad \forall t \in (0, T^*).
    \end{equation*}
    \\
    \\
    (\textit{C}) 
    By point (\textit{A}), $\abs{\ubi \spc}^2$ is a sub-solution in $\overline{\Real^3 \setminus \Omega_t} \times (0,T^*)$.
    Consequently, the weak parabolic maximum principle applies, and for all \, $t \in (0,T^*)$, the points where $\abs{\ubi(x,t) \spc}^2$ attains its absolute maximum
    lie on \mbox{$\partial \spc (\Real^3 \setminus \Omega_t) \cap \partial \spc \Omega_t$} and absolute maximum is non-increasing with respect to $t$.
    This leads to
   \begin{equation}\label{BoundNSE04}
        \norm{\ubi(\cdot,t)}_{\infty, \overline{\Real^3 \setminus \Omega_t}} \spc  \leq \norm{\ubi(\cdot,\varepsilon)}_{\infty, \overline{\Real^3 \setminus \Omega_{\varepsilon}}} \quad \forall t \, \in [\varepsilon, T^*).
    \end{equation}
    Moreover, by (\textit{B}), we have
    \begin{equation}\label{BoundNSE06}
      \norm{u(\cdot,t)}_\infty \leq \norm{\ubi(\cdot,\varepsilon)}_\infty \quad \forall t \in [\varepsilon, T^*).
    \end{equation}
    Since \(\lim_{t \to T^*} \norm{\ubi(t)}_\infty < \infty\), by Theorem \ref{T:KATO_NSE} \textit{(ii)}, the solution \(\ubi\) can be extended for all \(t \geq 0\).
    \\
\end{proof}

\medskip
\medskip
\medskip
\medskip
\medskip
\medskip

\vspace{1em}
\noindent\textbf{arXiv Classes:}

\noindent Primary Classes: Mathematics (math)

\noindent Subclasses: math.AP, math-ph/math.MP, math.DS, math.NA, physics.flu-dyn

\end{document}